\documentclass{article}
\usepackage{amsthm}
\usepackage{amsfonts}
\usepackage{cite}
\usepackage{amsmath, amscd}
\newtheorem{theorem}{Theorem}[section]

\newtheorem{proposition}[theorem]{Proposition}
\newtheorem{lemma}[theorem]{Lemma}

\begin{document}
\author{Jeremy Berquist}
\title{A Semi-Smooth Kodaira Vanishing Theorem and Demi-Normal Cyclic Covering Lemma}
\maketitle

\noindent
\linebreak
\textbf{Abstract.  }   The classical Kodaira Vanishing Theorem states that $H^i(X, \omega_X \otimes \mathcal{L}) = 0$ for $i>0$, where $X$ is a smooth projective variety over $\mathbb{C}$ and $\mathcal{L}$ is an ample line bundle on $X$.  We prove an analogous vanishing result under the assumption that $X$ is a semi-smooth projective variety over $\mathbb{C}$.  This is done by passing to the normalization of $X$, which is a smooth projective variety when $X$ is semi-smooth.  The classical vanishing result is proved using a cyclic cover.  Therefore we include a proof that the cyclic cover of a demi-normal variety is again demi-normal.

\newpage
\tableofcontents
\addcontentsline{}{}{}
\newpage
\begin{section}{Introduction}  The Kodaira Vanishing Theorem is used in birational geometry in order to prove (among other things) the Cone Theorem of the Minimal Model Program (MMP).  See \cite{KM98} Section 2.4 for the statement of the theorem and discussion on its proof.   Vanishing results are often useful when studying global properties of algebraic varieties.  Within the context of the MMP, it becomes important to study singular varieties, even when the data that begin the program are smooth varieties.  In particular, it becomes necessary to work also with non-normal varieties. 

The non-normal varieties that arise in the context of the MMP are the demi-normal varieties.  These are varieties $X$ with Serre's $S_2$ property, such that the codimension one part has only double normal crossing singularities.  Equivalently, $X$ is $S_2$, Gorenstein in codimension one ($G_1$), and seminormal (SN).  See \cite{GT80}.  The appropriate way to study the singularities on $X$ is to use a semi-resolution of $X$, in the same way that a resolution of singularities is a tool to study the singularities of a normal variety.  

A semi-resolution $f: Y \rightarrow X$ is a proper, birational morphism, such that the only singularities on $Y$ are double normal crossings and pinch points.  (These are desribed algebraically as points whose local rings have completions isomorphic to $k[[x_1, \ldots, x_n]]/(x_1x_2)$ for double normal crossing points and $k[[x_1, \ldots, x_n]]/(x_1^2 - x_2^2x_3)$ for pinch points.  We also require that a semi-resolution be such that $f$ sends the singular locus of $Y$ birationally onto the singular locus of $X$ in an open set where the only singularities are double normal crossings.  In particular, no component of the singular locus of $Y$ is $f$-exceptional.)  A variety $Y$ with these singularities is called semi-smooth.

The normalization $p:  \overline{Y} \rightarrow Y$ of a semi-smooth variety $Y$ is smooth.  Moreover, the conductor loci, given by $C$ in $\overline{Y}$ and by $D$ in $Y$, are such that $C$ and $D$ are both smooth, and $C \rightarrow D$ is a smooth morphism of degree two, ramified along the pinch locus.  Whenever we have a finite morphism such as $p$, ample line bundles pull back to ample line bundles.  Moreover, there is a simple relationship between the dualizing sheaves on $Y$ and $\overline{Y}$ (note that $Y$ has only hypersurface singularities, and is therefore Gorenstein, so its dualizing sheaf is invertible):  $p^*\omega_Y = \omega_{\overline{Y}}(C).$

It seems natural then that a vanishing result should hold for $Y$ when a similar result holds for its normalization.  Our main result is the following (2.1):  Given a semi-smooth variety projective variety $X$ over $\mathbb{C}$, and an ample line bundle $\mathcal{L}$ on $X$, the higher cohomology $H^i(X, \omega_X \otimes \mathcal{L}) =0$ for $i>0$.  note that Serre duality holds for $X$, so an equivalent statement is that $H^i(X, \mathcal{L}^{-1}) = 0$ for $i < \textnormal{dim } X$.

The proof of the Kodaira Vanishing Theorem uses cyclic covers.  Although we do not duplicate the proof here by substituting ``semi-smooth" for ``smooth" at the appropriate places (a method that is not sure to work), but instead pass to the normalization and use the statement of the theorem for smooth varieties, it seems appropriate to discuss cyclic covers in the demi-normal case.  This we do after proving the main result.  In particular, we show that a cyclic cover of a demi-normal variety is again demi-normal.  In principle, such a result allows us to, for example, replace a $\mathbb{Q}$-Gorenstein (Cohen-Macaulay) demi-normal variety with its index-1 cover, and thus to assume the Gorenstein condition.  We have yet to find a use for such a trick in birational geometry, because the cover thus constructed is not birational to the variety we start with.  However, its existence does suggest that Kodaira vanishing in the semi-smooth case should be a natural consequence of vanishing in the smooth case.

\end{section}

\begin{section}{Main Result}  Let $X$ be a semi-smooth projective variety over $\mathbb{C}$.  Let $p:  \overline{X} \rightarrow X$ be the normalization, so that $\overline{X}$ is smooth.  Note that the normalization of a projective variety is projective.  If the conductor loci are denoted by $C \hookrightarrow \overline{X}$ and $D \hookrightarrow X$, then the induced morphism $C \rightarrow D$ is a (smooth and finite) morphism between smooth varieties, of generic degree two.  In particular, the trace map provides a splitting $$\mathcal{O}_D \rightarrow p_*\mathcal{O}_C \rightarrow \mathcal{O}_D.$$  We use this fact by applying cohomology functors to $\mathcal{O}_C$ and $\mathcal{O}_D$ in what follows.

We have the following short exact sequences: (*)  $$   0 \rightarrow \mathcal{O}_X(-D) \rightarrow \mathcal{O}_X \rightarrow \mathcal{O}_D \rightarrow 0$$ and similarly for the ideal sheaf of $C$, noting that $p$ is finite, (**)$$    0 \rightarrow p_*\mathcal{O}_{\overline{X}}(-C) \rightarrow p_*\mathcal{O}_{\overline{X}} \rightarrow p_*\mathcal{O}_C \rightarrow 0.$$  Now the conductor ideal sheaf is by definition the largest such ideal sheaf in $\mathcal{O}_X$ that is already an ideal sheaf in $\mathcal{O}_{\overline{X}}$.  Hence the first terms in these sequences are isomorphic.  Moreover, there is a commutative diagram of exact sequences sending the first into the second.

Now let $\mathcal{L}$ be an ample line bundle on $X$.  Then $p^*\mathcal{L}$ is ample.  We have by duality that $p^*\omega_X = \omega_{\overline{X}}(C).$  We apply Kodaira vanishing.  It states that $H^i(\overline{X}, p^*\omega_X(-C) \otimes p^*\mathcal{L}) = 0$ for $i >0$.  Since $p$ is finite, there are no nonzero higher direct images, and the Leray spectral sequence shows that $$0 = H^i(X, p_*(p^*\omega_X(-C) \otimes p^*\mathcal{L})) = H^i(X, \omega_X \otimes \mathcal{L} \otimes p_*\mathcal{O}_{\overline{X}}(-C))$$ for $i >0$, where we've used the projection formula in the second equality.

In other words, when we tensor the sequence (**) by $\omega_X \otimes \mathcal{L}$ and take the long exact sequence in cohomology, we have isomorphisms for all $i > 0$:  $$H^i(X, p_*\mathcal{O}_{\overline{X}} \otimes \omega_X \otimes \mathcal{L}) \cong H^i(X, p_*\mathcal{O}_C \otimes \omega_X \otimes \mathcal{L}).$$

But the first term in (*) is isomophic to the first term in (**), so we conclude that we also have isomorphisms in cohomology from tensoring (*) with $\omega_X \otimes \mathcal{L}.$  Since there is a morphism of short exact sequences, as mentioned above, we obtain the following commutative diagram for all $i>0$:
$$\begin{CD}
H^i(X, \omega_X \otimes \mathcal{L})                                                                @>\cong>>               H^i(X, \omega_X \otimes \mathcal{L} \otimes \mathcal{O}_D) \\
@VVV                                                                                                                                                      @VVV \\
H^i(X, p_*\mathcal{O}_{\overline{X}} \otimes \omega_X \otimes \mathcal{L})  @>\cong>>             H^i(p_*\mathcal{O}_C \otimes \omega_X \otimes \mathcal{L}) \\
\end{CD}$$
The second vertical arrow splits, since there are maps $\mathcal{O}_D \rightarrow p_*\mathcal{O}_C \rightarrow \mathcal{O}_D$ whose composition is the identity.  Tensoring and applying $H^i$ are functorial, so we obtain the splitting as required.   In particular, the first vertical arrow also has a splitting map, meaning that the following composition is the identity:  $$H^i(X, \omega_X \otimes \mathcal{L}) \rightarrow H^i(X, p_*\mathcal{O}_{\overline{X}} \otimes \omega_X \otimes \mathcal{L}) \rightarrow H^i(X, \omega_X \otimes \mathcal{L}).$$  Thus it suffices to show that the middle term is zero for $i>0$.

We have $$H^i(X, p_*\mathcal{O}_{\overline{X}} \otimes \omega_X \otimes \mathcal{L}) = H^i(X, p_*(p^*\omega_X \otimes p^*\mathcal{L})),$$ by the projection formula, which since $p$ is finite equals $$H^i(\overline{X}, p^*\omega_X \otimes p^*\mathcal{L}) = H^i(\overline{X}, \mathcal{O}_{\overline{X}}(C) \otimes \omega_{\overline{X}} \otimes p^*\mathcal{L}),$$ where we've again used duality.  Thus we are done if we can show that this last term is zero.

Consider the exact sequence $$0 \rightarrow \mathcal{O}_{\overline{X}} \rightarrow \mathcal{O}_{\overline{X}}(C) \rightarrow \mathcal{O}_C(C) \rightarrow 0.$$  If we tensor by $\omega_{\overline{X}}$, adjunction implies that the final nonzero term is $\omega_C$.  Then when we tensor by $p^*\mathcal{L}$, we have $$0 \rightarrow \omega_{\overline{X}} \otimes p^*\mathcal{L} \rightarrow \mathcal{O}_{\overline{X}}(C) \otimes \omega_{\overline{X}} \otimes p^*\mathcal{L} \rightarrow \omega_C \otimes p^*\mathcal{L} \rightarrow 0.$$  We want to show that the higher cohomology of the middle term vanishes.  But this follows from the long exact sequence in cohomology and Kodaira vanishing for the smooth varieties $\overline{X}$ and $C$.  

We have therefore proved:

\begin{theorem}  Let $X$ be a projective semi-smooth variety over $\mathbb{C}$.  Then $$H^i(X, \omega_X \otimes \mathcal{L}) = 0$$ for $i>0$ and any ample line bundle $\mathcal{L}$ on $X$.

\end{theorem}

\end{section}
\begin{section}{Cyclic Covering Lemma}  In this section we show that a cyclic cover of a demi-normal variety is again demi-normal.  The relevance for the present situation is that Kodaira vanishing uses the cyclic covering trick.  Since a similar vanishing theorem holds in the non-normal situation, it seems appropriate that such a covering should be demi-normal.

\begin{lemma}  The property that a finite morphism is \'etale is local on the base.  In other words, a finite morphism $f: Y \rightarrow X$ is \'etale if and only if there is a covering of $X$ by affine open sets $U$ such that $f|_{f^{-1}U}: f^{-1}U \rightarrow U$ is \'etale.
\proof  By definition, $f$ is \'etale if and only if $f$ is flat and $\Omega_{Y/X} = 0$.  By construction of the relative cotangent bundle, we have $\Omega_{Y/X}|_{f^{-1}U} \cong \Omega_{f^{-1}U/U}$ for affine open sets $U$.  Note that $f^{-1}U$ is also affine in this case.  That $f$ is flat if and only if $f|_{f^{-1}U}$ is flat for all $U$ in an open covering of $X$ follows directly from the definition.  The lemma follows from these two observations.\qed
\end{lemma}

\begin{proposition}  Let $X$ be a $\mathbb{Q}$-Gorenstein variety with properties $G_1, S_2$, and seminormality.  There exists a finite cover $\pi: Y \rightarrow X$ in which $Y$ has the same properties, the morphism $\pi$ is \'etale over the locus where $K_X$ is Cartier, and $K_Y = \pi^*K_X$ is a Cartier divisor on $Y$.  If $Y$ is Cohen-Macaulay, then it is Gorenstein.
\proof  Let $m$ be an integer such that $mK_X$ is Cartier.  For each open set where $(\omega_X^{\otimes m})^{\vee\vee}$ (corresponding to $mK_X$) is trivial, we choose a generating section $s : \mathcal{O}_X \rightarrow \mathcal{O}_X(mK_X)$.  Then we construct the sheaf of $\mathcal{O}_X$-algebras $$\mathcal{A} = \mathcal{O}_X \oplus \mathcal{O}_X(K_X) \oplus \cdots \oplus \mathcal{O}_X((m-1)K_X),$$ with multiplication defined by $s$:  $$\mathcal{O}_X(aK_X) \otimes \mathcal{O}_X(bK_X) \rightarrow \mathcal{O}_X((a+b)K_X)$$ if $a+b <m$, or $$\mathcal{O}_X(aK_X) \otimes \mathcal(bK_X) \rightarrow \mathcal{O}_X((a+b)K_X) \stackrel{s^{-1}}{\rightarrow} \mathcal{O}_X((a+b-m)K_X)$$ if $a+b \geq m$.  Then consider the natural morphism $\pi: Y = \textbf{\textnormal{Spec}}_X \mathcal{A} \rightarrow X$.  

If this is done for each set in an open cover on which $mK_X$ is trivial, then we claim that the relative spectra and their projections patch together to give a finite morphism onto $X$.  The required properties of $Y$ are all local, and being \'etale is local on the base by (3.1).  Thus we may assume that $s: \mathcal{O}_X \rightarrow \mathcal{O}_X(mK_X)$ is a nowhere vanishing section giving a trivialization of $mK_X$.

To prove the claim, recall that we are working over an algebraically closed field of characteristic zero.  Then the construction above is independent of the choice of generating section, up to isomorphism.

The divisor $\pi^*K_X$ is Cartier.  Suppose without loss of generality that $K_X$ is effective (multiplying by an invertible sheaf does not change invertibility properties).  Then the inclusion maps $\mathcal{O}_X((i-1)K_X) \hookrightarrow \mathcal{O}_X(iK_X)$, together with the isomorphism $s^{-1} : \mathcal{O}_X(mK_X) \rightarrow \mathcal{O}_X$, defines an $\mathcal{A}$-linear map $g: \mathcal{A} \rightarrow \mathcal{A}$ such that $g^m$ is the local equation of $\pi^*(mK_X)$.  (An $R$-linear map from a ring $R$ to itself is completely determined by where it sends 1.  In this case, the identity is $(1,0,\ldots,0)$.  Multiplication in $R$ corresponds to composition of maps.)  Then $g$ can be thought of as a section of $\mathcal{O}_Y$, and $\pi^*K_X$ is the Cartier divisor given by $g$.  In other words, the rational function $g$ defines an invertible subsheaf of $\mathcal{K}_Y$, generated by the same element everywhere.  Thus $\pi^*K_X$, defined as $\frac{1}{m}\pi^*(mK_X)$, is equal to the divisor determined by $g$.

If $\mathcal{O}_X(K_X)$ is invertible, then locally on $X$ it is generated by a single rational function.  Thus over an affine chart $\textnormal{Spec}(B)$, $Y$ is given by the spectrum of $C = B[y]/(y^m-s)$, where (by slight abuse of notation) $s \in B$ is a unit.  The extension $B \rightarrow C$ is flat since it is free, and \'etale because the cotangent bundle can be seen to be zero.  For this, consider the exact sequence $$I/I^2 \stackrel{\delta}{\rightarrow} \Omega_{B[y]/B} \otimes_{B[y]} C \rightarrow \Omega_{C/B} \rightarrow 0.$$  Let $s'$ be an inverse of $s$ in $B$.  Then we have $$\delta(\frac{s'y}{m} \cdot (y^m - s)) = \frac{s'}{m}(y^m -s)dy \otimes 1 + \frac{s'y}{m} \cdot my^{m-1} dy \otimes 1.$$  Since we can move everything in the tensors but the $dy$ to the right, we see that we get $dy \otimes s'y^m \equiv dy \otimes 1$ in the image of $\delta$.  This element generates the second term in the above sequence, so by exactness we conclude that $\Omega_{C/B} = 0$.  Again since \'etaleness is local, $\pi$ is \'etale over the open set where $K_X$ is Cartier.

In particular, $\pi$ is \'etale over an open set $U \subset X$ whose complement has codimension at least two.  Since $\pi$ is finite, the preimage $\pi^{-1}U$ also has a complement with codimension at least two.  It follows that $K_Y$ and $\pi^*K_X$ agree in codimension one (since the relative cotangent bundle is zero there), and hence must be equal.  Thus $K_Y$ is a Cartier divisor.

The morphism $\pi$ is such that $\pi_*\mathcal{O}_Y = \mathcal{A}$.  The summands of $\mathcal{A}$ are $S_2$, being divisorial sheaves, hence so is their sum.  Then $Y$ is also $S_2$.  We claim that at every point $y \in \pi^{-1}U$, the local ring $\mathcal{O}_{Y,y}$ is Gorenstein and seminormal.  In fact, we may choose $U$ so that all the local rings at points of $U$ are Gorenstein and seminormal.  For any \'etale morphism, the completions of corresponding local rings differ only by a separable field extension:  $\hat{\mathcal{O}}_x \otimes_{k(x)} k(y) \cong \hat{\mathcal{O}}_y$.  A local ring is Gorenstein if and only if so is its completion.  See \cite{BH98} 3.1.19.  By \cite{GT80} 5.3, the same statement is true when the property is seminormality.  Finally, tensoring by a separable field extension of the residue field does not change either one of these properties.  For seminormality, this is shown in \cite{GT80} 5.7.  For the Gorenstein property, this follows from the Ext definition, along with the fact that change of fields defines an exact functor.  

Thus $Y$ has $G_1$ and $S_2$ and is seminormal in codimension one.  By \cite{GT80} 2.7, the $S_2$ condition implies that $Y$ is seminormal everywhere.

We note that being Gorenstein is the same as having an invertible dualizing sheaf and being Cohen-Macaulay.  We can always conclude that $K_Y$ is invertible, but we need to assume that $Y$ is Cohen-Macaulay before concluding that it is Gorenstein.  If all of the sheaves in $\pi_*\mathcal{O}_Y$ are Cohen-Macaulay, then so is $Y$.  See \cite{KM98}, 5.4.  However, even if $X$ itself is Cohen-Macaulay, it is not clear that the summands in $\pi_*\mathcal{O}_Y$ will be Cohen-Macaulay.  It seems the best result in this direction that we can prove is that $X$ is Cohen-Macaulay if $Y$ is.  This follows from the splitting $\mathcal{O}_X \rightarrow \pi_*\mathcal{O}_Y \rightarrow \mathcal{O}_X$ and the fact that Cohen-Macaulayness can be checked by the vanishing of local cohomology modules.\qed
\end{proposition}

We note that at least one other vanishing result for normal varieties and their resolutions holds for demi-normal varieties and their semi-resolutions.  Namely, Grauert-Riemenschneider vanishing.  As proved in \cite{Berq14}, for a semi-resolution $f: Y \rightarrow X$ of a demi-normal variety $X$, the higher direct images $R^if_*\omega_Y$ are zero.  This result is related to Kodaira vanishing in the normal case, via \cite{KM98} 2.68.
\end{section}
\nocite{BH98}
\nocite{GT80}
\nocite{KM98}
\nocite{Berq14}

\bibliographystyle{plain}
\bibliography{paper}

\end{document}